\documentstyle[12pt]{article}

\begin{document}

\newcommand{\ai} {\vbox to 7pt{\hbox to 7pt{\vrule height 7pt width 7pt}}}

\newcommand{\complement}{\cal C}

\newcommand{\hr} {\hookrightarrow}

\newcommand{\rt} {\rightarrow}

\newcommand{\st} {\stackrel}

\newcommand{\dps} {\displaystyle}

\newcommand{\pn} {\par\noindent}

\newcommand{\C} {\mbox{$I\kern-0.60emC$}}

\newcommand{\Z} {\mbox{$Z\kern-0.40emZ$}}

\newcommand{\lo} {\longrightarrow}

\newcommand{\sai} {\mbox{$\to \kern -0.50 em \to$}}

\newcommand{\nsai} {\mbox{$\not\to \kern -0.50 em \to$}}

\newcommand{\N} {{\rm I}\!{\rm N}}

\newcommand{\R} {{\rm I}\!{\rm R}}

\newcommand{\F} {{\rm I}\!{\rm F}}

\newcommand{\hs} {\hskip}

\newcommand{\vs} {\vskip}

\newcommand{\cqd} {\hspace{10pt} \rule {5pt}{5pt}}

\newcommand{\call} {{\cal C}}

\newcommand{\ov} {\overline}

\baselineskip=20pt

\vs 1cm

\begin{center}

\pn {\large {\bf  Some results about the Schroeder-Bernstein Property for separable Banach spaces}}

\vs 0.2cm

\normalsize

\pn {\small {by}}

\vs 0.2cm

\pn \small {Valentin Ferenczi and El\'oi Medina Galego }

\end{center}

\normalsize

\vs 0.3cm

\small {{\bf Abstract.} We construct a continuum of mutually non-isomorphic
  separable  Banach spaces which are complemented in each other. 
 Consequently, the Schroeder-Bernstein Index of any of these spaces is
$2^{\aleph_0}$. Our
  construction is based on a Banach space introduced by W. T. Gowers and
  B. Maurey in 1997.
 We also  use classical descriptive set theory methods, as in some work of
  V. Ferenczi and C. Rosendal, to improve  some results of P. G. Casazza and
  of N. J. Kalton on the
 Schroeder-Bernstein Property for
  spaces with  an unconditional finite-dimensional Schauder decomposition.}

\vs 0.5cm

\normalsize

\pn {\bf 1. Introduction}

\vs 0.2cm

 Let $X$ and $Y$ be Banach spaces. We  write $X \st{c}{\hr} Y$ if $X$ is
 isomorphic to a complemented subspace of $Y$, $X \sim Y$ if $X$ is isomorphic
 to $Y$ and $X \not \sim Y$ when $X$ is not isomorphic to $Y$. We also write
 $X \st{c} {\sim} Y$  if both $X \st{c}{\hr} Y$ and $Y \st{c}{\hr} X$ hold.
If $n \in \N  = \{1, 2, 3, \cdots \}$, then $X^n$  denotes the sum of $n$
 copies of $X$.  The first infinite cardinal number will be indicated by
 $\aleph_{0}$. We shall write our proofs in the case of real Banach spaces,
 clearly the results hold in the complex setting as well.

The Schroeder-Bernstein Problem for Banach spaces asks whether
isomorphism and complemented biembeddability must coincide for any pair of Banach
spaces. In other words, if $X$ and $Y$ are Banach spaces such that  $X \st{c} {\sim} Y$, does it follow that $X \sim Y$?
It was answered by the negative by Gowers \cite{G}, using separable
spaces. Later on, Gowers and Maurey \cite{GM2}
provided other counterexamples to the Problem: in particular, they built a separable Banach
space $X_1$, which is
isomorphic to its cube but not to its square. So $X_1 \st{c} {\sim} X_1^2$ while
$X_1 \not\sim X_1^{2}$.

The answer by the negative given by Gowers to the Schroeder-Bernstein Problem
opens two directions of research which are the guidelines of this paper.

\vs 0.3cm

{\bf 1.1. The Schroeder-Bernstein Property for Banach spaces with unconditional
 Schauder decomposition.}

\vs 0.2cm

 The first direction of research is to ask what additional conditions ensure a positive answer to the
Schroeder-Bernstein Problem. More precisely, a  Banach space $X$ is said to have the Schroeder-Bernstein Property (in short,
the SBP), if whenever a Banach space $Y$ satisfies $Y \st{c} {\sim} X$, it follows that
$Y \sim X$.

 We wish to find
sufficient conditions on a Banach space to have the SBP.
For example, according to the well known Pelczynski's Decomposition Method, a Banach space of
the form $l_p(X)$ for $1 \leq p <+\infty$ has the SBP.
For more information  about Banach spaces having the SBP, see the survey of
P. Casazza \cite{C}, and for more examples of Banach spaces failing the SBP, see \cite{MG}. 

We recall that a Banach space $X$ is said to have a
 {\em Schauder decomposition} $X=\sum_{n \in \N}E_n$, where 
$(E_n)_{n \in \N}$ is a sequence of
 closed subspaces of $X$, if every $x \in X$ can be written in a unique
way as $x=\sum_{n \in \N} x_n$, with $x_n \in E_n$ for all $n$.
It is {\em unconditional} if there exists a constant $C$ such that for all
$x=\sum_{n \in \N} x_n$, and every subset $I$ of $\N$, we have that
$\|\sum_{n \in I} x_n\| \leq C\|x\|$. 
A {\em finite-dimensional Schauder decomposition (or FDD)} is a Schauder decomposition
$X=\sum_{n \in \N}E_n$ for which $E_n$ is finite-dimensional for all $n$.
We shall use the classical abbreviation {\em UFDD} for an unconditional finite-dimensional Schauder decomposition.

Two important open problems about the $SBP$ are to know whether every
primary Banach space or every Banach space with an unconditional basis has the
$SBP$. Recall that a Banach space $X$ is said to be  {\em primary} if whenever $X=Y \oplus Z$, then 
$Y \sim X$ or $Z \sim X$.
 The following definitions introduced by  N. J. Kalton \cite{K}
are the starting point of our research in this  direction.

A Banach space $X$ has  the {\it SBP restricted to spaces with UFDD},
 if $X$ has an UFDD and whenever a Banach space $Y$ with UFDD satisfies $Y \st{c} {\sim} X$, it follows
that $Y \sim X$. A Banach space $X$ is said to be {\it countably primary} if there is a countable set
$(S_n)_{n \in \omega}$ of Banach spaces such that whenever $X= A \oplus B$,
then there exists $n \in \omega$ such that
 $A \sim S_n$ or $B \sim S_n$.

Kalton obtained various results about  
countably primary Banach spaces with an unconditional basis or
an unconditional Schauder decomposition.
In particular, a countably primary Banach space with an unconditional basis
has the SBP restricted to spaces with UFDD, see  Theorem 2.8 in \cite{K}.

In \cite{FR}, most results of Kalton concerning Banach spaces with an
unconditional basis were improved. The method was simplified, using classical
results of descriptive set theory. Uniformity of the constants of isomorphism
was obtained.
The results were also extended to
 $\kappa$-primary for any $\kappa<2^{\omega}$ (with the
obvious definition).

In Section 2, after noting
that a primary Banach space with a UFDD must have the SBP
(Proposition 2.1),  we show how to extend the methods of \cite{FR} to the case of
spaces with an unconditional Schauder decomposition. In consequence, we improve in a
similar way the
results of Kalton about spaces with an unconditional Schauder decomposition, and
our method is also more direct. One
application related to SBP is that for $\kappa<2^{\omega}$,
a $\kappa$-primary Banach space with an
unconditional basis has the SBP restricted to spaces with UFDD (Theorem 2.10).
 In fact, it is more generally true of Banach spaces with an unconditional
 basis which are not {\em perfectly decomposable} (roughly speaking, a Banach
 space $X$ is perfectly decomposable if there are perfectly many mutually
 non-isomorphic ways of decomposing $X$). We shall define precisely this
topological notion in Section 2.

The new uniformity result we get allows us to
improve on some work of Casazza concerning primary spaces: we prove
 that an $l_p$-sum of
finite-dimensional spaces, which is
$\kappa$-primary for some $\kappa<2^{\omega}$ (or more generally, which is  not perfectly decomposable)
  must have the
SBP (Corollary 2.12).

\vs 0.3cm

{\bf 1.2. The Schroeder-Bernstein Index for Banach spaces.}

 \vs 0.2cm

 The second
direction of research about the SBP is the following. Given the existence of non trivial
families of Banach spaces which are mutually non isomorphic yet
complementably biembeddable in each other, we wish to know what are the possible structures
for these families, for example in terms of cardinality. This question was formalized by
the definition of the Schroeder-Bernstein index $SBi(X)$ of a Banach space
$X$
in \cite{E}.
Here, we shall use a modified Schroeder-Bernstein Index
$SBI(X)$ given by Definition 1.1 below. This definition is simpler and more
natural than  the one of $SBi(X)$. 
Both indexes are cardinal numbers, and denoting  by $\alpha^+$  the successor cardinal of any cardinal $\alpha$, it is direct to check that they are
related by
 $SBi(X)=SBI(X)^+$.

\vs 0.2cm

{\bf Definition 1.1.} Let $X$ be a Banach space. Let $CB(X)$ be the set of
subspaces
$Y$ of $X$ such that $Y \st {c} {\sim} X$. Let
$\overline{CB}(X)=CB(X)/\sim$ be the set of isomorphism
 classes of elements of $CB(X)$.
Then
\vspace{-0,2cm}
$$SBI(X):=|\overline{CB}(X)|. \vspace{-0,2cm}
$$

\vs 0.1cm

More simply said, $SBI(X)$ is the number of $\sim$-classes on the
$\st {c} {\sim}$-class of $X$.

\vs 0.1cm

 Observe that for any Banach space
 $X$, $SBI(X)=1$ if and only if $X$ has the SBP. Also,
$SBI(X) \leq 2^{dens(X)}$, where  $dens(X)$ denotes
the density character of $X$.
 Indeed there is a basis of open sets
 for the topology of $X$, of cardinality $dens(X)$. So there are no more than
$2^{dens(X)}$ open subsets of $X$, and in particular no more than
$2^{dens(X)}$
(isomorphism classes of) closed subspaces of $X$.
 
\vs 0.1cm

Our goal is to see what are the possible values of 
the Schroeder-Bernstein Index for a Banach space.

Clearly the counterexample $X_G$ of Gowers in \cite{G}                                            satisfies $SBI(X_G)>1$.
 In 1997, Gowers and Maurey (\cite{GM2}, page 559) constructed, for each $p
 \in \N$, $p \geq 2$, a Banach space $X_1({\cal S}_p)$ (this notation will be
 explicited in Section 4), which is isomorphic to its
 subspaces
of codimension $n$ if and only if $p$ divides $n$.  Consequently
 $SBI(X_1({\cal S}_p))
 \geq p$,
 consider the family of spaces $(X_1({\cal S}_p) \oplus {\R}^n)_{0 \leq n \leq p-1}$. 
In fact, using the properties of $X_1({\cal S}_p)$ mentioned in the remark in \cite{GM2} after Theorem 19, it is not difficult to prove
that $SBI(X_1({\cal S}_p))=p$.

Recently the second author \cite{E} found a Banach space $X_E$ such that $X_E^2$ is isomorphic to a complemented subspace of  $X_E$, but 
$X_E^m$ is not isomorphic to $X_E^n$, for every $m \neq n$. Hence $(X_E^n)_{n \in
 \N}$ is an infinite sequence of mutually non-isomorphic Banach spaces which
 are complemented in each other. Thus $SBI(X_E) \geq {\aleph_0}.$

The main aim of Section 3, and the main result of this paper, is to provide a family of cardinality the continuum
of  mutually non-isomorphic separable Banach spaces which are complemented in
each other (Theorem 3.9). In particular, if $X$ is any  member of
this family then $SBI(X)=2^{\aleph_0}$. The construction of  such a  family is
inspired by the construction by Gowers and Maurey of the Banach space
$X_1$ isomorphic to its cube but not to its square  (\cite{GM2}, Section 4.4).

In Section 4, we note some open problems about the SBP as well as some side
consequences of some of our techniques. For example, we show that failing the
SBP is not a three-space property (Proposition 4.3).

 Finally, we end with an appendix which contains the proof of two
technical lemmas needed in Section 3.

\vs 0.5cm

{\bf 2.  On the SBP restricted to  Banach spaces with unconditional
  finite-dimensional Schauder decomposition.}

\vs 0.2cm

We start by noting a direct and interesting consequence of Kalton's results.

\vs 0.3cm

{\bf Proposition  2.1.} {\it Let $X$ be a primary Banach space with a UFDD.
 Then $X$ has the $SBP$.}

\vs 0.3cm

{\bf Proof.}
Being primary, $X$ must be isomorphic to its hyperplanes.
By  Theorem 2.3 from \cite{K}, for some $N \in\N$,
$X \sim X \oplus \sum_{n \geq N} E_n$, where the $E_n$'s are the summands of
the UFDD of $X$. From this and from the fact that each $E_n$ is
finite-dimensional, we deduce that $X$ must
be isomorphic to its square.

But a primary space isomorphic to its square is easily seen to have the SBP.
Indeed,  assume $Y \st {c} {\sim} X$. There exists $Z$ such that $Y \sim X \oplus Z$,
so $Y \sim X \oplus Z \sim X \oplus X \oplus Z \sim X \oplus Y$.
On the other hand $Y$ embeds complementably in $X$, so by the primariness
 of $X$, either $Y \sim X$, and we are done, or $X \sim X \oplus Y \sim Y$.
\cqd 

\vs 0.3cm

We now turn to our generalization of Kalton's result about the SBP restricted
to spaces with UFDD, for countably primary Banach spaces. We shall use
classical results and definitions from descriptive set theory, and our
reference for these will be the book of Kechris \cite{Ke}.
\vs 0.3cm

Let $X$ be a separable Banach space with a Schauder
decomposition  $\sum_{n=1}^{+\infty} E_n$.
 We start as in
 \cite{K} or \cite{FR} by assigning to each element $\alpha$ of $2^{\omega}$
 a
subspace $X(\alpha)$ of $X$ in the obvious way:
$$X(\alpha)=\sum_{\alpha(n)=1} E_n.$$

 The relation of isomorphism between spaces of the form $X(\alpha)$, 
$\alpha \in 2^{\omega}$, 
induces a relation on $2^{\omega}$ that we shall denote by $\simeq$, and it is
not difficult to check that it is analytic.

 We first
give some definitions.
For a cardinal $\kappa \leq 2^{\omega}$, we
 say that a Schauder decomposition is {\em $\kappa$-primary} if there is 
a set $\cal S$ of Banach spaces, of cardinality $\kappa$, such that for every
subset $I$ of $\N$, there exists $S \in {\cal S}$ such that
$\sum_{n \in I}  E_n \sim S$
or $\sum_{n \in \N \setminus I}  E_n \sim S$. 
A Banach space $X$ is {\em $\kappa$-primary} if 
if there is 
a set $\cal S$ of Banach spaces, of cardinality $\kappa$, such that 
whenever $X \sim A \oplus B$, there exists $S \in {\cal S}$ such that
$A \sim S$
or $B \sim S$.
Evidently an unconditional Schauder decomposition of a
$\kappa$-primary Banach space is $\kappa$-primary.

 \vs 0.1cm

We shall say that a Schauder decomposition is {\em  perfect}
if there is a perfect subset
 $P$ of $2^{\omega}$ such that,
 for all $\alpha, \beta \in P$ with $\alpha \neq \beta$,
$$\sum_{\alpha(n)=1}  E_n \not\sim \sum_{\beta(n)=1}  E_n,$$
$$\sum_{\alpha(n)=0}  E_n \not\sim \sum_{\beta(n)=1}  E_n,$$
$$\sum_{\alpha(n)=0}  E_n \not\sim \sum_{\beta(n)=0}  E_n.$$
\hs 0.6cm In particular, note that a Schauder decomposition which is
 $\kappa$-primary for some 
$\kappa < 2^{\omega}$, is not perfect.

For our last definition, we need to recall that the set of separable Banach
spaces, seen as subspaces of an isometrically universal separable Banach space
such as $C([0,1])$, or more generally, the set of subspaces of a given
separable Banach space $X$, may be equipped naturally with a Borel structure called
the Effros-Borel structure (see e.g. \cite{FR}). In this setting we may talk
about Borel or analytic sets of separable Banach spaces, or of subspaces of a
given separable Banach space $X$; note that any
uncountable Borel set of Banach spaces is necessarily of cardinality
$2^{\omega}$.

For a Banach space $X$, we call {\em decomposition of $X$}  a pair $(A_0,A_1)$
of subspaces
of $X$ such that $X=A_0 \oplus A_1$. 
We say that a separable Banach space $X$ is {\em  perfectly decomposable} if
 there is a Borel set
 $\{(A_{\alpha}^0$, $A_{\alpha}^1), \alpha \in 2^{\omega}\}$
of decompositions of $X$, such that for $\alpha \neq \beta$ and
any $(\epsilon,\gamma)$ in $\{0,1\}^2$,
$A_{\alpha}^{\epsilon} \not\sim A_{\beta}^{\gamma}$.
So a Banach space   which is
 $\kappa$-primary for some 
$\kappa < 2^{\omega}$, is  not perfectly decomposable.

Finally and evidently, if a separable Banach space has a perfect
 unconditional Schauder decomposition then it is  perfectly decomposable.

\vs 0.3cm

We need to recall two theorems from descriptive set theory
(Theorems 19.1 and 8.41 from \cite{Ke}).

\vs 0.3cm

{\bf Theorem 2.2 (Kuratowski-Mycielski)} {\it Let $E$ be a perfect Polish space,
and $R$ be a relation on $E$ which is meager in $X^{2}$. Then there exists a
homeomorphic copy $C$ of the Cantor space such that $\forall x,y \in C$ with
$x \neq y$, we have $x \neg R y$.}

\vs 0.3cm

{\bf Theorem 2.3 (Kuratowski-Ulam)} {\it Let $E$ be a Polish space and $D$ be a
subset of $E^2$ having the Baire property. Then
$D$ is nonmeager if and only if
$$\exists^* x \exists^* y: (x,y) \in D.$$}

Here $\exists^*x P(x)$ signifies the existence of a nonmeager set of $x$ such
that $P(x)$.

\vs 0.2cm

We are now ready to prove a proposition in the spirit of \cite{FR}.

\vs 0.3cm

{\bf Proposition 2.4.} {\it Let $X$ be a separable 
Banach space with a Schauder
decomposition $X=\sum_{n=1}^{+\infty}  E_n$ which is not perfect. Then
there exists an $\simeq$-class which is non-meager in
  $2^{\omega}$.}

\vs 0.2cm

{\bf Proof.} For $\alpha \in 2^{\omega}$, we shall denote by $\complement \alpha$
the element $(1-\alpha(n))_{n \in \omega}$ of $2^{\omega}$.
As in \cite{FR} we define the relations $\simeq_1$ and $\simeq_2$ on $2^{\omega}$ by
$$\alpha \simeq_1 \beta \Leftrightarrow \alpha \simeq \complement \beta,$$
$$\alpha \simeq_2 \beta \Leftrightarrow \complement \alpha \simeq \complement \beta.$$

The first case in our proof is to assume that $\simeq$, $\simeq_1$ and $\simeq_2$
are meager. Then their union is meager as well. We then apply Theorem 2.2. to
get a perfect set $P$ avoiding this union, i.e. with the property stated in the definition of perfect Schauder decompositions.

In the second case, assume
for example that $\simeq_2$ is non-meager. Being analytic, $\simeq_2$ has the
Baire property, so by Theorem 2.3., we may find an element
$\alpha \in 2^ {\omega}$ such that, for $\beta$ in a non-meager subset of
$2^{\omega}$, $\complement \alpha \simeq \complement \beta$. As clearly, the
map
sending  $\beta$ to $\complement \beta$ is an homeomorphism on $2^{\omega}$, we deduce that the
$\simeq$-class of $\complement \alpha$ is non-meager.
A similar proof holds if $\simeq$ or $\simeq_1$ is non-meager. $\cqd$

\vs 0.3cm

{\bf Remark 2.5.} In \cite{FR}, it was shown that in the case of Banach space
with a $1$-dimensional Schauder decomposition (i.e. with a Schauder basis),
a non-meager $\simeq$
class in $2^{\omega}$ must be comeager. This used the fact that
modifying a finite number of vectors of the basis
of a Banach space preserves the isomorphism
class. We cannot use this fact in the general case of a Schauder
decomposition, and in fact the result is  false in that case:
consider $X=l_1 \oplus (\sum l_2)_{l_2}$.
We see that the $\simeq$-class corresponding to $l_2$ is non-meager
 but not comeager.

\

We now deduce from Proposition 2.4 the following extension of
Theorem 3.4 from \cite{K}. For $X$ and $Y$ Banach spaces, and
$K \geq 1$, $X \sim^K Y$ means that $X$ is $K$-isomorphic to $Y$.

\vs 0.3cm

{\bf Theorem 2.6.} {\it Let $X$ be a  separable Banach space with an unconditional Schauder
decomposition $X=\sum_{n=1}^{+\infty} E_n$ which is not perfect. Then,
there exists an integer $N$ and a constant $K$
such that, for every subset $I$ of
$[N,\infty)$, $X \sim^K X \oplus (\sum_{n \in I} E_n.)$}

\vs 0.2cm

{\bf Proof.} By Proposition 2.4, we assume there is a non-meager $\simeq$-class
 and we intend to find $K$ and $N$ such that for every subset $I$ of
$[N,\infty)$, $X  \sim^K X \oplus \sum_{n \in I} E_n$.
Let $X(\alpha_0)$ for some $\alpha_0$ be a Banach space
in the isomorphism class associated to the non-meager $\simeq$-class.
For $n \in \N$, let ${\cal A}_n$ be the
set of $\alpha$'s such that $X(\alpha) \sim^n X(\alpha_0)$.
Then for some $K \in \N$, $\cal A={\cal A}_K$ is non-meager.

The set $\cal A$ is analytic, thus has the Baire property, and we deduce that
$\cal A$ is comeager in some basic open set $N(u)$. Here $u=(u(1),\ldots,u(k))$ denotes an
element of $2^{<\omega}$. In other words the set $\cal C$ of $\alpha \in 2^{\omega}$
such that the concatenation $u^{\frown}\alpha \in \cal A$ is comeager in $2^{\omega}$.
 We now apply the following  classical characterization of comeager subsets of
$2^{\omega}$
which was already used (and proved) in \cite{FR}. As $\cal C$ is comeager in $2^{\omega}$,
then there exists a partition of
$\N$ in infinite subsets $M_1$ and $M_2$, subsets $N_1 \subset M_1$
and $N_2 \subset M_2$, such that: 
for $i=1,2$,
an element $\alpha$ of $2^{\omega}$ is in $\cal C$ whenever
for every $n \in M_i$, $\alpha(n)=1$ if $n \in N_i$ and
$\alpha(n)=0$ if $n \not\in N_i$.

 Going back to $\cal A$, we get a partition of
$[k+1,+\infty)$ in infinite subsets $A_1$ and $A_2$, subsets $B_1 \subset A_1$
and $B_2 \subset A_2$, such that: 
for $i=1,2$,
an element $\alpha$ of $2^{\omega}$ is in $\cal A$ whenever

(a) for every $n=1,\ldots,k$, $\alpha(i)=u(i)$, and

(b) for every $n \in A_i$, $\alpha(n)=1$ if $n \in B_i$ and
$\alpha(n)=0$ if $n \not\in B_i$.

\vs 0.1cm

Now let $N=k+1$, let $I$ be any subset of $[N,\infty)$, and
let $Y=\sum_{i \in I}\oplus E_i$.

We let $C_1= A_1 \cap I$,
$C_2= A_2 \cap I$.
The element $\alpha_1$ defined by $\alpha_1(i)=u(i)$, for $i<N$, and
$\alpha_1(i)=1$ if and only if $i \in C_1 \cup B_2$, for $i \geq N$, satisfies
(a) and (b), so
belongs to $\cal A$. Denoting by $X(u)$ the finite sum $\sum_{u(i)=1} E_i$,
we deduce that
\vspace{-0,2cm}
 $$X(\alpha_0) \sim X(u) \oplus (\sum_{i \in C_1 \cup B_2} E_i). \vspace{-0,2cm}
$$
\hs 0.6cm  By unconditionality of the basis,
\vspace{-0,2cm}
$$X(\alpha_0) \sim X(u) \oplus  (\sum_{i \in B_2} E_i)
\oplus (\sum_{i \in C_1} E_i). \vspace{-0,2cm}
$$
\hs 0.6cm By the characterization of $\cal A$ again,
\vspace{-0,2cm}
$$X(\alpha_0) \sim X(\alpha_0) \oplus (\sum_{i \in C_1} E_i). \vspace{-0,2cm}
$$
\hs 0.6cm Likewise
$$X(\alpha_0) \sim X(\alpha_0) \oplus (\sum_{i \in C_2} E_i).
 $$
\hs 0.6cm Combining the two, and noting  that $C_1$, $C_2$ form a partition of $I$,
we get
\vspace{-0,3cm}
 $$X(\alpha_0) \sim X(\alpha_0) \oplus Y. \vspace{-0,3cm}$$
\hs 0.6cm  As $X(\alpha_0)$ is complemented in $X$, we
deduce
\vspace{-0,3cm}
$$X \sim X \oplus Y. \vspace{-0,1cm}$$
\hs 0.6cm Finally, a look at the proof shows that we can get this isomorphism with an uniform
constant depending on $K$ and the constant of unconditionality of the Schauder
decomposition. \cqd

\vs 0.3cm

Concerning Theorem 2.6,  note that a stronger conclusion  of the form $X \sim X^2$ (as in the
case of a space with an unconditional basis) is false: consider a decomposition
$X=X_0 \oplus (\sum l_2)_{l_2}$, where $X_0$ is the
 hereditarily indecomposable Banach space of Gowers and Maurey \cite{GM1}.
There are only four classes of isomorphism of subspaces $X(\alpha)$
(including the classes of $\{0\}$ and $X_0$), but it
is not difficult to show that $X$ is not isomorphic to its square.

However, as an immediate consequence of Theorem 2.6, we have:

\vs 0.3cm

{\bf Corollary 2.7.} {\it Suppose that the  separable Banach space $X$ has 
Schauder unconditional
  decomposition which is not perfect, and assume also that  every summand 
is isomorphic to its square. Then $X$ is isomorphic to its square.}

\vs 0.3cm

The next theorem is a variation of Theorem 2.6.
We recall that $E_0$ is the equivalence
 relation defined on $2^{\omega}$ by $\alpha E_0
\beta$ if and only if $\exists m: \forall n \geq m, \alpha(n)=\beta(n)$. It is
the $\leq_B$-lowest Borel equivalence relation above equality on $2^{\omega}$.
For the definition of $\leq_B$: given an equivalence relation $R$ on a Polish
space
$E$, (resp. $R'$ on $E'$), we say that 
$(E,R)$ is Borel reducible to $(E',R')$, and write
$(E,R) \leq_B (E',R')$, if there is a Borel map $f:E \rightarrow E'$, such
that
for all $x,y$ in $E$, $x R y$ if and only if $f(x) R' f(y)$. 
We refer to \cite{FR2} for more about the notion of the relation $\leq_B$ of Borel reducibility
between equivalence relations.
When a relation reduces $E_0$, in particular there is a perfect set of mutually
non-related points.

In \cite{FR2}, a Banach space $X$ was defined to be {\em ergodic} if the relation $E_0$
is Borel reducible to isomorphism between subspaces of $X$. It was proved by
Ferenczi and Rosendal (\cite{FR},\cite{R}) that a Banach space $X$ with an
unconditional basis which is not ergodic, must be isomorphic to its square, to
its hyperplanes and more generally to $X \oplus Y$ for any subspace $Y$
generated by a finite or infinite subsequence of the basis.

We can prove a result in the same spirit for Banach spaces with a UFDD,
provided we (unessentially)
relax the ergodic assumption.  
We recall that Banach spaces $X$ and $Y$ are {\em nearly isomorphic}, and
write  $X \st{f} {\sim} Y$, if some finite-codimensional subspace of $X$ is
isomorphic to a finite-codimensional subspace of $Y$. We shall say that a Banach space is {\em  nearly ergodic} if $E_0$ is Borel
reducible to $\st{f} {\sim}$ between subspaces of $X$.
 The next proposition shows that when a Banach
space
is ergodic, then it is nearly ergodic. Note    that both imply that there is a perfect
set of mutually non (nearly) isomorphic subspaces. 

\vs 0.3cm

{\bf Proposition 2.8.} {\it Let $X$ be an ergodic Banach space. Then it is nearly ergodic.}

\vs 0.2cm

{\bf Proof.} By definition there exists a Borel map $g$ from $2^{\omega}$ into
the set of subspaces of $X$, such that $\alpha E_0 \beta$ if and only if
$g(\alpha) \sim^f g(\beta)$. Denote by $E_f$ the relation defined on
$2^{\omega}$ by $\alpha E_f \beta$ if and only if
$g(\alpha) \st{c} {\sim} g(\beta)$. An $E_f$-class is a countable union of $E_0$-classes
and thus is meager. As every $E_f$-class is meager, then by Theorem 2.3, 
$E_f$ is meager in $(2^{\omega})^2$. 
Every $E_f$-class is also invariant by $E_0$, so by Proposition 14 of
\cite{R}, the relation $E_0$ is Borel reducible to $E_f$. Combining the map
reducing $E_0$ to $E_f$ with
$g$, we get a Borel reduction of $E_0$ to $\st{f} {\sim}$ on the set of subspaces of
$X$. \cqd

\vs 0.3cm

{\bf Theorem 2.9.} {\it Let $X$ be a separable Banach space with an unconditional Schauder
decomposition $X=\sum_{n=1}^{+\infty} E_n$ such that $E_n$ is of finite dimension for $n>N$, for some $N$. Assume
$X$ is not nearly ergodic.
Then there exists an integer $k$ and a constant $K$ such that, for every subset $I$ of
$[k,\infty)$, $X \sim^K X \oplus (\sum_{n \in I} E_n.)$}

\vs 0.2cm

{\bf Proof.} Let  $\st{f} {\simeq}$ be the relation induced on
$2^{\omega}$ by near isomorphism between spaces of the form $X(\alpha)$,
$\alpha \in 2^{\omega}$. 

 First
 assume $\st{f} {\simeq}$ is meager. Let $1_N$ be the length $N$ sequence
$(1,\ldots,1)$. Then the relation $r$, defined on $2^{\omega}$ by
\vspace{-0,3cm}
$$\alpha r \beta \Leftrightarrow 1_N^{\frown}\alpha \st{f} {\simeq}
1_N^{\frown}\beta, \vspace{-0,3cm}$$
is also meager in $(2^{\omega})^2$. Furthermore, because $E_n$ is of finite
dimension for $n>N$, $r$ is clearly invariant by a
finite change of coordinates of the sequences $\alpha$ and $\beta$.
 According to Proposition 14 of \cite{R}, we
deduce that $E_0$ is Borel reducible (by some $g$) to $r$. 
The map $f$ from $2^{\omega}$ into $2^{\omega}$, defined by
$f(\alpha)=1_N^\frown g(\alpha)$, is then a Borel reduction of $E_0$
to $\st{f} {\simeq}$, that is, $X$ is nearly ergodic, a contradiction.

So assume $\st{f} {\simeq}$ is non meager. Using Theorem 2.3 as before, we deduce
that some $\st{f} {\simeq}$ class is non-meager.
This class is a countable union of $\simeq$-classes, so
we deduce that some $\simeq$ class is non-meager. We may then proceed
as in Theorem 2.6.  \cqd

\vs 0.3cm

 The result of Kalton about the SBP restricted to spaces with UFDD, for
countably primary Banach spaces, mentioned in the introduction, now
generalizes
 to Banach spaces which are either not perfectly decomposable  or not nearly ergodic.

\vs 0.3cm

{\bf Theorem 2.10.} {\it Let $X$ be a Banach space with an unconditional
  basis. Assume $X$ is not perfectly decomposable or not nearly ergodic.
Then $X$ has the SBP restricted to spaces with UFDD.}

\vs 0.2cm

{\bf Proof.} By Theorem 2.6 or Theorem 2.9, there exists $N \in \N$ such that 
for every $I \subset [N,+\infty)$, $X \sim X \oplus \sum_{i \in I} \R e_i$.
Taking $I=\{N\}$ we see that $X$ is isomorphic to its hyperplanes.
Taking $I=[N,+\infty)$ it follows that $X$ is isomorphic to its square. 
 
Now let $Y$ be a Banach space   with an UFDD $(E_n)_{n \in \N}$
such that $Y \st{c} {\sim} X$ .
Then there exists a space $W$ such that $Y \sim X \oplus W$ and we deduce
\vspace{-0,3cm}
$$Y \sim X \oplus W \sim X \oplus X \oplus W \sim X \oplus Y. \vspace{-0,3cm}$$ \hs 0.6cm Also
there exists some space $Z$ such that
\vspace{-0,3cm}
$$X \sim Z \oplus Y \sim Z \oplus (\sum_{n \in \N} E_n). \vspace{-0,2cm}
$$ 
\hs 0.6cm This UFDD satisfies the hypotheses of Theorem 2.6 or Theorem 2.9. We deduce
that for some $K \in \N$, $X \sim X \oplus \sum_{n \geq K}E_n$.
As the decomposition is finite-dimensional and
$X$ is isomorphic to its hyperplanes, it follows that
\vspace{-0,3cm}
$$X \sim X \oplus Y. \vspace{-0,3cm} \cqd$$

\vs 0.3cm

In \cite{Ca}, Casazza proved that if a Banach space $X$ is an $l_p$-sum of 
finite-dimensional spaces and is primary, then $X \sim l_p(X)$, and thus $X$
has the SBP. The following Corollary 2.12 extend this result to
$l_p$-sums of finite-dimensional spaces which are not perfectly decomposable
or not ergodic. We point out here that to prove Corollary 2.12, we need uniformity
in the result of Proposition 2.11. This uniformity is one of the results we
get which was not proved in the paper of Kalton.

\vs 0.3cm

{\bf Proposition 2.11.} {\it Let $X$ be a Banach space with a UFDD
$X=\sum_{n=1}^{+\infty} E_n$. Assume this UFDD is not perfect or
the relation $E_0$ is not Borel reducible to near isomorphism between
subspaces of $X$ of the form
 $X(\alpha)=\sum_{\alpha(n)=1}E_n, \alpha \in 2^{\omega}$ (for example $X$
 could be
non nearly ergodic). Then there exists $N \in \N$ and an infinite sequence of
disjoint subsets $(B_k)_{k \in \N}$ of $[N,+\infty)$  such that if 
$Y=\sum_{n=N}^{+\infty} E_n$, then $Y \sim \sum_{k \in \N} (\sum_{i \in B_k} E_i)$, with $\sum_{i \in B_k} E_i$  uniformly isomorphic to Y.}

\vs 0.2cm

{\bf Proof.} By Proposition 2.4 or the proof of Theorem 2.9, we know that
some
$\simeq$-class is non-meager. We start as in the proof of Theorem 2.6, and
using the same notation. Let $\alpha_0$ be in some fixed non-meager
$\simeq$-class, we may find $K$ such that the set $\cal A$ of $\alpha$'s such
that $X(\alpha) \sim^K X(\alpha_0)$ is non-meager, and thus comeager in some
$N(u)$, $u \in 2^{\omega}$.
We note that, by adding a finite sum of spaces $E_i$, and up to modifying the 
constant $K$ and $\alpha_0$, we may assume that $u(i)=1$ for all $i \leq |u|$
and also that $\alpha_0 \in N(u)$.

As a preliminary result, let us prove that the non-meager $\simeq$-class we
get corresponds to the isomorphism class of $X$. We let $X(u)=\sum_{i \leq |u|} E_i$.
Let $N=|u|+1$. We apply the same proof as in Theorem 2.6 to get, for all
$Y=\sum_{i \in I} E_i$, $I \subset [N,+\infty)$:

$$X(\alpha_0) \sim X(\alpha_0) \oplus Y.$$
\hs 0.6cm In particular,
$$X(\alpha_0) \sim X(\alpha_0) \oplus \sum_{i \geq N}E_i,$$
and thus $$X(u) \oplus X(\alpha_0) \sim X(\alpha_0) \oplus X.$$
\hs 0.6cm On the other hand, we also have for all $Y=\sum_{i \in I} E_i$, $I \subset
[N,+\infty)$:
$$X \sim X \oplus Y.$$
\hs 0.6cm Choose $Y$ such that $X(\alpha_0)=X(u) \oplus Y$, then
$$X \oplus X(u) \sim X \oplus X(\alpha_0).$$

Finally, we deduce that $X(u) \oplus X(\alpha_0) \sim X(u) \oplus X$, and as
$X(u)$ is finite-dimensional, that $X(\alpha_0) \sim X$.

So, modifying $K$ if necessary, we may assume that $\alpha_0(i)=1,\ \forall i$,
that is $X(\alpha_0)=X$. Now we prove the result about the decomposition. 
 We note that the characterization of comeager subsets of $2^{\omega}$ in 
terms of partitions of $\N$, that we used in Theorem 2.6, can be generalized
 to an infinite partition (see \cite{FR} about this).
 So from the fact that the set of $\alpha$'s such that $X \sim^K X(\alpha)$
is comeager in $N(u)$, we get a sequence $(B_n)_{n \in
 \N}$
of disjoint
 subsets of $[N,+\infty)$ such that:

\vs 0.1cm

(a) $X \sim^K X(u) \oplus \sum_{i \in B_k} E_i$, for all $k \in \N$,
and

(b) $X \sim X(u) \oplus \sum_{k \in \N} (\sum_{i \in B_k} E_i)$.

\vs 0.1cm

Let $Y=\sum_{n=N}^{+\infty} E_n$. From (b), we have that
$$X(u) \oplus Y=X \sim X(u) \oplus \sum_{k \in \N} (\sum_{i
  \in B_k} E_i),$$
and thus as $X(u)$ is finite-dimensional,
$$Y \sim \sum_{k \in \N} (\sum_{i \in B_k} E_i).$$ 
\hs 0.6cm From (a), we get
$$X(u) \oplus Y= X \sim^K X(u) \oplus (\sum_{i \in B_k} E_i),
\forall k \in \N,$$
and so,
$$Y \sim (\sum_{i \in B_k} E_i), \forall k \in \N,$$
where the constant of  isomorphism depends only on
$K$ and on $\dim X(u)$. $\cqd$

\vs 0.3cm

{\bf Corollary 2.12.} {\it Let $Z=c_0$ or $l_p$ for $1\leq p<+\infty$.
 Let $X=(\sum_{n \in \N}E_n)_Z$, where for each $n$, 
$E_n$ is finite-dimensional. Assume this UFDD is not perfect or $X$ is not ergodic. Then $X \sim (\sum X)_Z$. So by Pelczynski's Decomposition Method, $X$ has the SBP.}

\vs 0.2cm

{\bf Proof.} Let $X$ be as in the hypotheses. First note that any $l_p$-sum, and in particular $X$ or any
subspace
of $X$ of the form $X(\alpha), \alpha \in 2^{\omega}$, contains a
complemented copy of $l_p$ and so, is isomorphic to its hyperplanes. 
In particular, two subspaces $X(\alpha)$ and $X(\beta)$,
 $\alpha,\beta \in 2^{\omega}$, are nearly isomorphic if and only if they are
isomorphic, and so the relation $E_0$ is Borel reducible to near isomorphism between
these subspaces if and only if it is Borel reducible to isomorphism between
them. This implies that the hypotheses of Proposition 2.11 are satisfied.

So  the conclusion of Proposition 2.11 holds,
and we note that the copies are uniform and
that the infinite direct sum is, in this case, an $l_p$-sum. So for some $N
\in \N$, $Y=\sum_{n \geq N} E_n$ satisfies $Y \sim l_p(Y)$.
As $X$ is isomorphic to its hyperplanes, it follows that $X \sim Y$, and so
$X \sim l_p(X)$. $\cqd$

\vs 0.5cm

{\bf 3. A continuum of mutually non-isomorphic Banach spaces which are complemented in each other.}

\vs 0.2cm

 We now turn to the main theorem of this paper, which provides a family of
cardinality the continuum  of mutually non-isomorphic separable
Banach spaces which are complemented in each other (Theorem 3.9). In order to present this family of Banach spaces, we need to fix  some notation and background from \cite{GM2}.

\vs 0.3cm

{\bf 3.1. Preliminaries.}

\vs 0.2cm

 Let $c_{00}$ be the vector space of all complex sequences which are eventually
$0$. Let $(e_n)_{n \in \N}$ be the standard basis of $c_{00}$. Given a vector $a=\sum a_n e_n$ its {\it support}, denoted supp(a), is the set of $n$ such that $a_n \neq 0$. Given subsets $E$, $F$ of $\N$, we say that $E<F$ if every element of $E$ is less than every element of $F$. If $x$, $y \in c_{00}$, we say that $x<y$ if $\hbox{supp}(x)<\hbox{supp}(y)$. If $x_1<x_2< \cdots<x_n$, then we say that the vectors $x_1$, $x_2$, $\cdots$ $x_n$ are {\it successive}. An infinite sequence of sucessive non-zero vectors is also called a {\it block basis} and a subspace generated by a block basis is a {\it block subspace}.

Given a subset $E$ of $\N$ and a vector $a$ as above, we write $Ea$ for the vector $\sum_{n \in E} a_n e_n$. An interval of integers is a set of the form $\{ m, m+1, \cdots, n \}$, where $m, n \in \N$. The {\it range} of a vector $x$, written $ran(x)$, is the smallest interval containing supp(x).

The following set of functions was first defined by Schlumprecht in \cite{S} except for condition (vi) which was added in \cite{GM2}.
 $\cal F$ denotes the set of functions ${\it f}: [1, \infty) \rt [1, \infty)$ satisfying the following conditions.

(i) ${\it f}(1)=1$ and ${\it f}(x)<x$ for every $x>1$;

(ii) {\it f} is strictly increasing and tends to infinity;

(iii) $\hbox{lim}_{x \rt \infty} x^{-q} {\it f}(x)=0$ for every $q>0$;

(iv) the function $x/{\it f}(x)$ is concave and non-decreasing;

(v) ${\it f}(xy) \leq {\it f}(x){\it f}(y)$ for every $x, y \geq 1$;

(vi) the right derivate of ${\it f}$ at 1 is positive.

\vs 0.1cm

Let $\cal X$ stand  for the set of normed spaces $(c_{00}, \| . \|)$ such that
the sequence $(e_n)_{n \in \N}$ is a normalized bimonotone basis, this means
that $\|Ex\| \leq \|x\|$ for every vector $x \in c_{00}$ and every interval
$E$. Given $X \in \cal X$ and ${\it f} \in \cal F$, it is said that $X$ {\it satisfies a lower f-estimate} if, given any vector $x \in X$ and any sequence of intervals $E_1 < E_2 < \cdots < E_n$, we  have  $\|x\| \geq {{\it f}(n)}^{-1} \sum_{i=1}^{n} \|E_{i}x\|.$

Given two infinite subsets $A$ and $B$ of $\N$, Gowers and Maurey define the {\it spread from A to B} to be the map $S_{A,B}:c_{00} \rt c_{00}$ defined as follows. Let $\rho   $ be the order-preserving bijection from  $A$ to $B$, then $S_{A,B}(e_n)=e_{{\rho}(n)}$ when $n \in \N$, and $S_{A,B}(e_n)=0$ otherwise.

Given any set  $\cal S$  of spreads, they say that it is a {\it proper set} it it is closed under composition and taking adjoints and if for every $(i,j) \neq (k,l)$, there are only finitely many spreads $S \in \cal S$ for which $e_{i}^{*}(Se_j) \neq 0$ and   $e_{k}^{*}(Se_l) \neq 0$.

\vs 0.1cm

For every Banach space $X$ satisfying a lower {\it f}-estimate for some ${\it
  f} \in \cal F$, and for every subspace $Y$ of $X$ generated by a block
basis, they define a seminorm $|||.|||$ on the set  $L(Y,X)$ of linear
mappings from $Y$ to $X$ as follows.  For $X \in \cal X$, and every integer $N
  \geq 1$, consider the equivalent norm on $X$ defined by
\vspace{-0,3cm}
$$\|x\|_{(N)}= \hbox{sup} \sum_{i=1}^{N} \|E_{i}x\|,\vspace{-0,3cm}
$$
where the supremum is extended to all sequences $E_1, E_2, \cdots, E_N$ of successive intervals.  Let ${\cal L}(Y)$ be the set of sequences $(x_n)_{n \in \N}$ of successive vectors in $Y$ such that $\|x_n\|_{(n)} \leq 1$,  for every $n \in \N$. Now let $$|||T|||=\sup_{x \in {\cal L}(Y)} \hbox{lim sup}_{n} \|T(x_n)\|.$$

Finally, we also  recall that two Banach spaces $X$ and $Y$ are said to be
totally incomparable if no infinite dimensional subspace of $X$  is  isomorphic to a subspace of $Y$.

\vs 0.3cm

{\bf 3.2. The main result.} 

\vs 0.2cm

The principal purpose of this section is to prove Theorem 3.9. Before that,
 we need
some auxiliary results which are similar to the ones involved in 
 the construction by
Gowers and Maurey of the Banach space $X_1$ isomorphic to its cube but not to its
square.  We shall improve these results in two directions.
First we shall need to have a family of spaces $X_r$ constructed on the model
of $X_1$, for $1/2 \leq r \leq 1$, and we shall take care that for $r \neq r'$,
$X_r$ and $X_{r'}$ are totally incomparable spaces. For this we shall have to write
new versions of some technical lemmas in \cite{GM2}.
Then, we shall need to know that each $X_r$ is not  nearly isomorphic to its
square (instead of just non-isomorphic). This requires a
little bit of extra care in the proofs as well.  

\vs 0.2cm

The proof of Lemma 3.1 is implicit in \cite{GM2}, page 864, and in the proof of \cite{GM2},
Lemma 9, in the case  $f(x)=\hbox{log}_{2}(x+1)$.

\vs 0.3cm

{\bf Lemma 3.1.} {\it Let  $f \in \cal F$ and  $J \subset \N$ be such
 that, if   $m, n \in \N$, $m<n$,   then
 $\hbox{log} \ \hbox{log} \ \hbox{log} \ n \geq 4m^{2}$. Write $J$ in
 increasing order as  $\{j_1, j_2, \cdots \}$ and let $K= \{j_1, j_3, j_5, \cdots \}$. Suppose that}

(a) {\it $f^{1/2} \in \cal F$;}

(b) {\it $f(j_1) > 256$;}

(c) {\it $\hbox{exp} \ \hbox{exp} \ j_{n} < f^{-1}(f(j_m))^{1/2}, \ \forall m, n \in \N, m<n$;}

(d) {\it $16 f(x^{1/2}) \geq f(x), \ \forall x \geq 0$;}

(e) {\it $f(j_n)^{3/2}>f(j_{n}^3), \ \forall n \in \N$;}

(f) {\it For every  $N \in J \setminus K$ and  $x_{0}$ in the interval} [log
N, exp N], {\it the function given by the tangent to $(x \rightarrow x/ f(x))$ at $x_0$ is at least $x/ f^{1/2}(x)$  for all positive $x$ outside the interval} $[\hbox{log} \ \hbox{log} \ N,  \ \hbox{exp} \ \hbox{exp} \  N]$.

\vs 0.1cm

 {\it Then for every  $K_0 \subset K$, there is a function $g \in \cal F$ such that $f \geq g \geq {f}^{1/2}$, $g(k)={f(k)}^{1/2}$ whenever $k \in K_{0}$ and $g(x)=f(x)$ whenever $N \in J \setminus K_{0}$ and $x$ in the interval} [log N, exp N].

\vs 0.3cm

The next  lemmas extend some technical results of \cite{GM2}
concerning
$f(x)=\log_2(x+1)$ to the case of $f(x)=(log_2(x+1))^{r}$, $ r \in [1/2,1]$.
Their calculus proofs are postponed until Section 5.
 
\vs 0.3cm

{\bf Lemma 3.2.} {\it Let $f_r$ be defined on $[1,+\infty)$ by $f_{r}(x)=(log_2 (x+1))^r$  for every $r \in (0, 1]$.
Then $f_r$ belongs to the class $\cal F$.}

\vs 0.3cm

{\bf Lemma 3.3.} {\it There exists  $J \subset \N$ such that writing it in increasing order as $\{j_1, j_2, \cdots \}$ and letting $K= \{j_1, j_3, j_5, \cdots \}$, we have that  for every   $K_0 \subset K=(j_{2i-1})_{i \in \N}$ and  every $r \in [1/2, 1]$, there is a function $g_r \in \cal F$ such that $f_r \geq g_r \geq {f_{r}}^{1/2}$, $g_{r}(k)={f_{r}(k)}^{1/2}$ whenever $k \in K_{0}$ and $g_{r}(x)=f_{r}(x)$ whenever $N \in J \setminus K_{0}$ and $x$ in the interval} [log N, exp N].

\vs 0.3cm

The case $r=1$ in the following theorem is the main result of [8], see [8, Theorem 5].  By  using Lemma 3.3 instead of [6, Lemma 6] in the argumentation of [8] we obtain:

\vs 0.3cm

{\bf Theorem 3.4.}     {\it Let $\cal S$ be a proper set of spreads and  $r \in [1/2, 1]$. There exists a Banach space $X_{r}(\cal S)$ satisfying a lower $f_{r}$-estimate with the following three properties.}

(i) {\it For every $x \in X_{r}(\cal S)$ and every $S_{A,B} \in \cal S$, $\|S_{A,B}x\| \leq \|x\|$ and therefore $\|S_{A,B}x\| = \|x\|$ if $\hbox{supp(x)} \subset A$.}

(ii) {\it If  $Y$ is a subspace of $X_{r}(\cal S)$ generated by a block basis, then every operator from $Y$ to $X_{r}(\cal S)$ is in the $|||.|||$-closure of the set of restrictions to $Y$ of the operators in the algebra $\cal A$ generated by $\cal S$. In particular, all operators on $X_{r}(\cal S)$ are $|||.|||$-perturbations of operators in $\cal A$.}

(iii) {\it The seminorm $|||.|||$ satisfies the algebra inequality $|||UV||| \leq  |||U||| \ |||V|||.$}

\vs 0.3cm

 We are now in position to define our family of totally incomparable 
versions of the
 'cube but not square' space $X_1$ of Gowers and Maurey.
 The definitions follow the ones in \cite{GM2}. For $i=0, 1, 2$ let  $A_i$ be
 the set of positive integers equal $i+1$ (mod 3), let $S'_i$ be the spread
 from $\N$ to $A_i$ and ${\cal S'}$ be the semigroup generated by $S'_0$,
 $S'_1$ and $S'_2$ and their adjoints. In \cite{GM2},
 page 560, it was shown that
 ${\cal S}'$ is a proper set.  Given $r \in [1/2, 1]$,  the Banach space  we
 are interested in is the space $X_{r}({\cal S}')$ obtained by Theorem 3.4.  As in \cite{GM2} it will be useful to   define it slightly less directly as follows.

Let $\cal T$ be the ternary tree $\cup_{n=0}^{\infty} \{0, 1, 2 \}^n$. Let $Y_{00}$ be the vector space of finitely supported scalar sequences indexed by $\cal T$ (including the empty sequence). Denote the canonical basis of $Y_{00}$ by $(e_t)_{t \in \cal T}$, write $e$ for $e_{\emptyset}$. If $s, t \in \cal T$, let $s^{\frown}t$ stand for the concatenation of s and t. We shall now describe some operators on $Y_{00}$.

Let $S_i$ for $i=0, 1,2$ be defined by their action on the basis as follows: $S_{i}e_t=e_{t^{\frown}i}$. The adjoint $S^*$ acts in the following way: $S_{i}^{*}e_t=e_s$ if t is of the form $t=s^{\frown}i$ and $S_{i}^{*}e_t=0$ otherwise.

If $P$ denotes the natural rank one projection on the line $\R e$, then we denote by  $\cal I$ and $\cal A$ respectively the proper set generated by $S_0$, $S_1$ and $S_1$, and the algebra generated by this proper set. Strictly speaking, $\cal I$ is not a proper set, but it is easy to embed $\cal I$ in $\N$ so that the maps $S_0$, $S_1$ and $S_1$ become spreads as defined earlier.

In order to obtain the space $X_{r}(\cal S')$, consider the subset ${\cal T}'$ of $\cal T$ consisting of all words $t \in \cal T$ that do not start with 0 (including the empty sequence). We modify the definition of $S_0$ slightly, by letting $S'_{0}e$ equal $e$ instead of $e_0$. The operators $S'_1$ and $S'_2$ are defined exactly as $S_1$ and $S_2$ were.

 To each $s=(i_1, \cdots i_n) \in {\cal T}'$ we can associate  the integer $n_s=3^{n-1}i_1+ \cdots +3{i_{n-1}}+i_{n}+1$, with $n_{\emptyset}=1$, and this defines a bijection between ${\cal T}'$ and $\N$. The operators $S'_{0}$, $S'_1$ and $S'_2$ coincide     with the spreads on $c_{00}$ defined earlier, so we can define ${\cal S}'$ to be the proper set they generate and obtain the space  $X_{r}({\cal S}')$.

 Let $Y$ be the completion of $Y_{00}$ equipped with the $l_1$ norm, in other words let $Y=l_{1}(\cal T)$ and let $\cal E$ denote the norm closure of $\cal A$ in $L(Y)$. Let also $\cal J$ the closed two-sided ideal in $\cal E$ generated by $P$.   Let $\cal O$ denote the quotient algebra ${\cal E}/ {\cal J}$.

Now we consider the algebra ${\cal A}'$ generated by ${\cal S}'$. In
 \cite{GM2}, Lemma 20, it was proved that $|||.|||$ is a norm on ${\cal
 A}'$. If we write $\cal G$ for the $|||.|||$-completion of ${\cal A}'$  then
 Theorem 3.4 implies that $\cal G$ is a Banach algebra and  there exists  a
 unital algebra homomorphism  $\phi$ from $L(X_{r}(\cal S'))$  to
 $\cal G$ (\cite{GM2}, page 550).

 In \cite{GM2}, Lemma 25, it was shown that there is  a norm-one algebra homomorphism $\theta$ from $\cal G$ to $\cal O$.

Finally, let $V$ be an arbitrary set and $\psi: V \rt V$ a function, we denote by ${\psi}_{3}$ the function from the set of matrices $M_{3}(V)$ to $M_{3}(V)$ given by ${\psi}_{3}((v_{i,j})_{1 \leq i,j \leq 3})=(\psi(v_{i,j}))_{1 \leq i,j \leq 3}$, for every $(v_{i,j})_{1 \leq i,j \leq 3} \in M_{3}(V)$.

 \vs 0.3cm
 
{\bf Remark 3.5.} Gowers and Maurey proved that $X_{1}({\cal S'})$ is
isomorphic to its cube ${X_{1}({\cal S'})}^3$ (\cite{GM2}, page 563). Likewise
$X_{r}({\cal S'})$ is isomorphic to its cube ${X_{r}({\cal
    S'})}^3$. Furthermore, it is
important to note that the norms of the projections involved in $X_{r}({\cal S'})     \st{c}{\hr} {X_{r}({\cal S'})}^2$ and ${X_{r}({\cal S'})}^2 \st{c}{\hr} X_{r}({\cal S'})$ do not depend on the number $r \in [1/2, 1]$. Indeed, letting
\vspace{-0,3cm}
$$X_{j}=\{\sum_{i=1}^{\infty} x_{3i+j} \  e_{3i+j} \in X_{r}({\cal S'}) \} \vspace{-0,3cm}$$
for every $j=0, 1, 2$, we have by Theorem 3.4 (i)

(a) $X_j$ is isometric to $X_{r}({\cal S'})$ for every $j=0, 1, 2.$

(b)$X_{r}({\cal S'})= X_0 \oplus X_1 \oplus  X_2.$

(c) The operator $P_r$ from $X_{r}({\cal S'})$ onto $X_0 \oplus  X_1$ defined by $P_{r}(x)=S_{\N, A_{0}}(x)+S_{\N, A_{1}}(x)$ is a projection with $\|P_r\| \leq 2$.

(d) The operator $Q_r$ from $X_0 \oplus  X_1$ onto $X_0$ defined by $Q_{r}=S_{\N,A_0}(x)$ is a projection with $\|Q_{r}\|=1$.

\vs 0.3cm

{\bf Lemma 3.6.} {\it $X_{r}({\cal S}')$ is  nearly isomorphic to ${X_{r}({\cal S}')}^2$  for no $r \in [1/2, 1]$.}

\vs 0.2cm

{\bf Proof.} It is inspired by \cite{GM2}, Theorem 26, where it was proved
 that  $X_{1}(\cal S')$ is not isomorphic
 to its square ${X_{1}({\cal S}')}^2$. Given $r \in [1/2, 1]$, denote  $X_{r}(\cal S')$ by $X$  and  assume that
 some finite codimensional subspace of  $X$ is isomorphic to a finite codimensional subspace of $X^2$. We consider the two possible cases.

\vs 0.1cm

{\it First case:}  $X \sim X^2 \oplus F$, for some finite-dimensional space $F$. 
 Let $U$ be an isomorphism from $X$ onto $X^2 \oplus F$, and assume without loss of generality that $F \subset
X$. Write $U=(U_1, U_2, U_3)$, where $U_1 \in L(X)$, $U_2 \in L(X)$ and $U_3 \in L(X, F)$.

Let $H$ be such that $X=F \oplus H$. Thus there exists an isomorphism $V$ from $X^2$ onto $H$. Defining $V_1(z)=V(x, 0)$ and $V_2(x)=V(0, x)$, for every $x \in X$, we have that $V_1 \in L(X)$,  $V_2 \in L(X)$ and $V(x_1, x_2)=V_{1}(x_1) + V_{2}(x_2)$, for every $(x_1, x_2) \in X^2$.

Next consider the isomorphism $\Psi:X \oplus X^2 \rt X^2 \oplus F \oplus H$ defined by
\vspace{-0,3cm}
$$\Psi(x_1, (x_2, x_3))=U(x_1) + V(x_2, x_3) \vspace{-0,3cm}$$
\hs 0.6cm Its matrix as linear map from $X^3$ to $X^3$ is given by
\vspace{-0,3cm}
$$A=\left(\begin{array}{lll}U_1 & 0 & 0\\
\vspace{-0,3cm} \\U_2 & 0 & 0\\
\vspace{-0,3cm} \\ U_3 & V_1 & V_2 \end{array}\right) \vspace{-0,3cm}$$
\hs 0.6cm where $U_3$ is seen as an operator from $X$ into $X$. Since   $A$ is an invertible element of $M_{3}(L(X))$, it follows that
${\theta}_3 {\Phi}_3 (A)$ is invertible in $M_{3}(\cal O)$. Therefore, by \cite{GM2}, Corollary 24, there exists an invertible element $B$ in $M_{3}(\cal E)$ such that ${\Pi}_{3}(B)= {\theta}_3 {\Phi}_3 (A)$, where  $\Pi$ is the canonical application from  $\cal O$ onto ${\cal E}/ {\cal J}$.

 As was noted in \cite{GM2}, page 550, the Kernel of $\phi$ is the set of $T \in
 L(X)$ satisfying $|||T|||=0$. The
 basis $(e_n)_{n \in \N}$  being shrinking (\cite{GM2}, page 551), the Kernel of $\phi$
 contains the compact operators. Indeed for any $x=(x_n)_{n \in \N} \in {\cal L}(Y)$,
 the
sequence $(x_n)_{n \in \N}$ is bounded and converges weakly to $0$. So if $T$ is
 compact,
$(T(x_n))_{n \in \N}$ converges in norm to $0$. In particular, $\phi(U_3)=0$, because  $U_3$ 
is of finite rank.

 On the other hand, according to \cite{GM2}, page  564,  $\cal J$ consists exactly of the compact $w^{*}$-continuous operators on $l_1$. Hence
\vspace{-0,3cm} 
$$B=\left(\begin{array}{ll} u & c_1\\
\vspace{-0,3cm} \\ c_2 & v \end{array}\right) \vspace{-0,3cm}$$
where $u \in M_{2,1}(\cal E)$, $v \in M_{1,2}(\cal E)$, $c_1 \in M_{1,1}(\cal
E)$ and $c_2 \in M_{2,2}(\cal E)$, with $c_1$ and $c_2$ compacts. It follows
from \cite{LT}, page 80, that the  operator $D$ below defined is Fredholm
$$D=\left(\begin{array}{ll} u & 0\\
\vspace{-0,3cm} \\ 0 & v \end{array}\right) \vspace{-0,3cm}$$

\vs 0.2cm

 Consequently $u$ and $v$ are also Fredholm operators; which is absurd because
 there exists no Fredholm element in $M_{2,1}(\cal E)$ (\cite{GM2}, Lemma 21).

\vs 0.1cm

{\it Second case:}  $X^2 \sim X \oplus F$  for some finite-dimensional space $F$.

In this case,  we  would have $X^3 \sim X^2 \oplus F$. According  to Remark 3.5,  $X \sim X^3$. Thus by the first case we also would obtain a contradiction.
$\cqd$

\

We recall the definition of {\em Rapidly Increasing Sequences} given in
\cite{GM2}. For $X \in \cal X$, $x \in X$ and every integer $N \geq 1$, recall that
\vspace{-0,3cm}
$$\|x\|_{(N)}= \hbox{sup} \sum_{i=1}^{N} \|E_{i}x\|,\vspace{-0,3cm}
$$
where the supremum is extended to all sequences $E_1, E_2, \cdots, E_N$ of successive intervals.

For $0< \epsilon \leq 1$ and $f \in \cal F$, we say that a sequence $x_1, x_2,
\cdots, x_N$ of successive vectors {\it satisfies the $RIS(\epsilon)$ condition
  for the function f} if there is a sequence
$(2N/f'(1))f^{-1}(N^{2}/{{\epsilon}^{2}})<n_1< \cdots <n_N$ of integers (where
f'(1) is the right derivate in $1$) such that $\|x_i \|_{(n_i)} \leq 1$ for each $i=1, \cdots, N$ and
\vspace{-0,3cm}
$$\epsilon {f(n_i)}^{1/2}>|ran (\sum_{j=1}^{i-1} x_j)| \vspace{-0,3cm}
$$ 
for every $i=2, \cdots, N$ (\cite{GM2}, page 546).

\vs 0.3cm

{\bf Lemma 3.7.} {\it Let ${\cal S}_1$ and ${\cal S}_2$ be proper sets of spreads and  $r, s \in [1/2, 1]$, with  $r \neq s$. Then the Banach spaces   $X_{r}({\cal S}_1)$ and $X_{s}({\cal S}_2)$ are totally incomparable.}

\vs 0.2cm

{\bf Proof.} Fix $r > s$ in $[1/2, 1]$  and suppose  that  $X_{r}({\cal S}_1)$ and
$X_{s}({\cal S}_2)$ are not totally incomparable. Thus, by a standard perturbation
argument, we may find  an  infinite sequence of successive non-zero vectors
$(z_n)_{n \in \N}$ in $X_r({\cal S}_1)$, and an isomorphism $T$ from
 $\overline {\hbox{span}} \{z_{n}: n \in \N\}$ into  $X_{s}({\cal S}_2)$, such that $(T(z_n))_{n \in \N}$ is successive in $X_{s}({\cal S}_2)$.

 For $N$ in $K \subset J$, we may then block $(z_n)_{n \in \N}$ to construct a
 sequence $x_1, x_2, \cdots, x_N$ satisfying the R.I.S(1)  condition for the
 function $f_r$ with $\|x_i\| \geq 1/2$ (\cite{GM2}, Lemma 4). Now putting
 $x=\sum_{n=1}^{N} x_n$ we obtain by the lemma analogous to Lemma 7 in \cite{GM2} that
\vspace{-0,3cm}
$$\|x\| \leq \frac{4 N}{f_r(N)}. \vspace{-0,3cm}$$ 
\hs 0.6cm Consequently 
$$\|T(x)\| \leq  \|{T}\| \frac{4 N}{f_r(N)}.$$
\hs 0.6cm On the other hand, since $T(x_n)$ is successive in $X_{s}({\cal S}_2)$ and   $X_{s}({\cal S}_2)$ satisfies a lower $f_{s}$-estimate we deduce that\vspace{-0,3cm}
$$\|{T(x)}\| \geq  \frac{1}{f_{s}(N)} {\sum_{n=1}^N \|T(x_n)\|} \geq   \frac{1}{f_{s}(N)} \frac{1}{\|{T^{-1}\|}} {\sum_{n=1}^N \|x_n\|}\geq \frac{N}{2}   \frac{1}{f_{s}(N)} \frac{1}{ \|{T^{-1}\|}}. \vspace{-0,3cm}$$
\hs 0.6cm It follows  that
 $$(\log_2 (N+1))^{r-s} \leq 8 \ \|{T}\| \ \|{T^{-1}\|},$$
which is a contradiction for $N$ large enough. $\cqd$

\vs 0.3cm

The last ingredient for our proof is the following standard lemma from \cite{HW}.

\vs 0.3cm

{\bf Lemma 3.8.} {\it Let $(Y_i)_{i \in \N}$ be a sequence of Banach
  spaces. Suppose that the Banach space $X$ is isomorphic to a subspace of
  $l_2(Y_i)_{i \in \N}$. Then some subspace of $X$ is isomorphic to a subspace of
  $Y_n$ for some $n \in \N$ or $l_2$ is isomorphic to a
  subspace of $X$.}

\vs 0.3cm

 We are now ready to prove our main result.

\vs 0.3cm

{\bf Theorem 3.9.} {\it There exists a family of separable Banach spaces  $(X_\alpha)_{\alpha \in {\{1,2\}}^{\omega}}$   such that}

(a) {\it $X_{\alpha} \st{c}{\hr} X_{\beta}$  for every $\alpha$ and $\beta$ in ${\{1,2\}}^{\omega}$.} 

(b)  {\it $X_{\alpha} \not \sim X_{\beta}$ for every $\alpha \neq \beta$.}

\vs 0.2cm

{\bf Proof.} Let $\cal S'$ be the spread considered after Theorem 3.4.  We pick a sequence  $(r_n)_{n \in \N}$ of  numbers in $[1/2,1]$,   with $r_m \neq r_n$ whenever $m \neq n$ and define   $Z_{n}=X_{r_n}(\cal S')$, for every ${n \in \N}$.

We then define for $\alpha \in {\{1,2\}}^\omega$, the following $l_{2}$-sum of Banach spaces:
\vspace{-0,3cm}
$$X_{\alpha}= (\sum_{n=1}^{\infty} Z_n^{\alpha(n)})_{2} \vspace{-0,3cm}.$$
\hs0.6cm It follows from Remark 3.5 that any two such spaces are complemented in each other.
We next assume that there exists  an isomorphism  $T$ from $X_{\beta}$ onto $X_{\alpha}$ and intend to prove that
$\alpha=\beta$. 
Let $P_1$ be the canonical projection from $X_\alpha$ onto $Z_1^{\alpha(1)}$,
and $Q_1$ be the canonical projection from $X_\alpha$ onto
$(\sum_{n=2}^{\infty} Z_n^{\alpha(n)})_{2}$.
We define analogously $P^{\prime}_1$ and $Q^{\prime}_1$ to be the projections corresponding to the space $X_\beta$.

We now claim that $S=Q_1 T_{| Z_1^{\beta(1)}}$ is strictly singular, that is, $S$ cannot be an isomorphism on any infinite dimensional subspace of $Z_1^{\beta(1)}.$

Indeed, otherwise, there exists an infinite dimensional subspace $Z$ of
$Z_1^{\beta(1)}$, such that $Q_1 T_{|Z}$ is an isomorphism into, so
$Z$ is isomorphic to some subspace of $(\sum_{n=2}^{\infty}
Z_n^{\alpha(n)})_{2}$. Since $Z_1^{\beta(1)}$ has no unconditional basic
sequence
 (\cite{GM2}, page 567), it contains no subspace isomorphic to $l_2$ and
 therefore   by Lemma 3.8, we 
deduce that $Z$ contains an infinite dimensional subspace which is isomorphic to a subspace of  $Z_n^{\alpha(n)}$ for some $n \geq 2$, contradicting  Lemma 3.7.

We then define the operators $U:Z_1^{\beta(1)} \rightarrow Z_1^{\alpha(1)}$ and  $V:Z_1^{\alpha(1)} \rightarrow Z_1^{\beta(1)}$ by
\vspace{-0,1cm} 
$$U(x)=P_1 T(x) \ \ \hbox{and} \ \ V(x)=P^{\prime}_1 T^{-1}(x). \vspace{-0,1cm}$$

 We consider $VU \in L(Z_1^{\beta(1)})$. For any $x \in Z_1^{\beta(1)}$, we deduce that
\vspace{-0,1cm} 
$$ VU(x)=VP_1 T(x)=V(Id_{X_{\alpha}}-Q_1)T(x)=V(T(x)-Q_1 T(x))=VT(x)-VS(x),\vspace{-0,1cm}$$
therefore
\vspace{-0,1cm}
$$VU(x)=P^{\prime}_1 (x) -VS(x)=
x-VS(x),\vspace{-0,1cm}$$
that is, $VU=Id_{Z_1^{\beta(1)}}+s$, where $s$ is strictly singular.

Now by symmetry, we also have that $UV=Id_{Z_1^{\alpha(1)}}+s'$, where $s'$ is
strictly singular. Then $UV$ and $VU$ are Fredholm operators (\cite{LT},
 page 80). It
follows that $U$ and $V$ are isomorphisms on finite codimensional subspaces
and have finite dimensional cokernels, hence are Fredholm. Consequently $Z_1^{\alpha(1)}$ and $Z_1^{\beta(1)}$ have isomorphic finite codimensional subspaces, and by Lemma 3.6, this means that $\alpha(1)=\beta(1)$. 
This proof can be repeated for an arbitrary $n \in \N$, so  we conclude that $\alpha=\beta$. $\cqd$

\vs 0.5cm

{\bf 4. Some remarks and problems.} 

\vs 0.2cm

{\bf Problem 4.1.} As already said, it comes easily from the
properties of the space $X_1({\cal S}_p)$ of Gowers
and Maurey mentioned in the introduction that it
satisfies $SBI(X_1({\cal S}_p))=p$.
Does there exist a Banach space $X$ with $SBI(X)=\aleph_0$? In particular, is the space $X_E$
defined in \cite{E} such a space?

\vs 0.3cm

{\bf Problem 4.2.} $SBI(X)=2^{\aleph_0}$ is the highest possible value
 for a separable Banach space $X$. The next step concerning separable Banach spaces
 should  rather be expressed
in terms of complexity of the relation of isomorphism.
We refer to \cite{FM} or \cite{R} for a survey about the notion of relative
complexity of analytic equivalence relations on Polish spaces, applied to isomorphism between
separable Banach spaces.  
How complex can be an equivalence relation $R$ on $2^{\omega}$, which is
 Borel reducible to isomorphism
between separable Banach spaces, with the condition that the image of the
 reducing map is formed by Banach spaces which are all complemented in each
 other?

\vs 0.3cm

As a consequence of Lemma 3.7, we also derive the following proposition.
 We recall that a property $P$ of a Banach space is said to be a
{\em three-space property} if whenever a Banach space $X$ has a subspace $Y$
which satisfies $P$ and such that $X/Y$ satisfies $P$, it follows that $X$
satisfies $P$. See \cite{CG} for a survey on three-space problems.

\vs 0.3cm

{\bf Proposition 4.3.} {\it Failing the Schroeder-Bernstein Property is not a three-space property.}

\vs 0.2cm

{\bf Proof.} Let $S$ be the right shift on $c_{00}$, that is, $S: c_{00} \rt c_{00}$ is given by $S(a_1, a_2, \cdots)= (0, a_1, a_2, \cdots)$. Denote by ${\cal S}_2$ and ${\cal S}_3$ the proper set generated by $S^2$ and $S^3$ respectively. Consider the Banach spaces $X=X_{1}({\cal S}_2)$ and  $Y=X_{1/2}({\cal S}_3)$ given by Theorem 3.4. By Lemma 3.7, $X$ and $Y$ are totally incomparable spaces.

We know by Theorem 19 in \cite{GM1}   and the remarks after this theorem that:

(a) $X$ is isomorphic to its subspaces of codimensions  two but not to its hyperplanes.

(b) Two finite-codimensional subspaces of $Y$ are isomorphic if and only if their codimensions are equal mod 3.

(c) Every complemented subspace of $X$ (resp. $Y$) has finite dimension or
codimension in $X$ (resp. $Y$).

 Clearly from (a) and (b) $X$ and $Y$ do not have the SBP. In fact $SBI(X)=2$ and $SBI(Y)=3$.  So, to prove the
 Proposition it suffices to show that $X \oplus Y$ has the SBP.

 Suppose then that $Z \st{c}{\hr} X \oplus Y$ and $X \oplus Y \st{c}{\hr}Z$
  for some Banach space $Z$. Since $X$ and $Y$ are totally incomparable
  spaces, by Theorem 23 in \cite{W}, we have that $Z=Z_1 \oplus Z_2$, where
  $Z_1 \st{c}{\hr} X$  and $Z_2 \st{c}{\hr} Y$. Moreover,  $X \st{c}{\hr} Z_1
  \oplus Z_2$  and $X$ and $Z_2$ are totally incomparable
  spaces. Consequently, again by Theorem 23 in \cite{W},  we conclude that
  $Z_1$ is an infinite dimensional space. According to (c) we deduce that $Z_1
  \sim X$ or $Z_1 \sim X \oplus \R$. In the same way, we obtain that $Z_2 \sim
  Y$, $Z_2 \sim Y \oplus \R$ or 
  $Z_2 \sim Y \oplus {\R}^2$.

Therefore $Z$  is isomorphic to one of the following  spaces: $X \oplus Y$, $X
\oplus Y  \oplus \R$, $X \oplus Y  \oplus {\R}^2$ or $X \oplus Y  \oplus
{\R}^3$. Hence, to see that $Z$ is isomorphic to $X \oplus Y$, it is enough to
show that $X \oplus Y$ is isomorphic to its hyperplanes.
 But this is true because of (a) and (b). Indeed,
\vspace{-0,3cm}
$$X \oplus Y \oplus \R \sim (X \oplus \R^2) \oplus Y \oplus \R \sim X \oplus
(Y \oplus \R^3) \sim X \oplus Y.  \vspace{-0,3cm} \cqd$$

\vs 0.3cm

Proposition 4.3 leads naturally to the following problems:
 
\vs 0.3cm

{\bf Problem 4.4.} Assume $X$ is a Banach space such that $X^2$ has the
SBP. Does it follow that $X$ has the SBP?

\vs 0.3cm

{\bf Problem 4.5.} Is the SBP a three-space property? 

\vs 0.3cm

A partial answer to Problem 4.5 was given by Casazza in \cite{C}. He noticed that if
$X$ and $Y$ have the SBP, and are totally incomparable spaces, then $X \oplus Y$ has
the SBP.
 
\vs 0.3cm

 As another immediate consequence of Theorem 3.4 and Lemma 3.7, 
we obtain the following result of I. Gasparis \cite{Ga}.  We recall that a Banach space $X$ is hereditarily indecomposable if no closed subspace $Y$ of $X$  contains a pair of infinite dimensional closed subspaces $M$ and $N$ such that $Y = M \oplus N$.

\vs 0.3cm

{\bf Corollary 4.6. (I. Gasparis)} {\it There exists a family of 
 cardinality the continuum of separable totally incomparable hereditarily indecomposable Banach spaces.}

\vs 0.2cm

{\bf Proof.} Let ${\cal S}=Id$ be the identity of $c_{00}$. Then, by
 \cite{GM2}, section 4.1, and  using Lemma 3.7 we see  that the spaces $(X_{r}(Id))_{r \in [1/2, 1]}$ given by Theorem 3.4 are a continuum of totally incomparable hereditarily indecomposable Banach spaces. $\cqd$

\vs 0.3cm

{\bf Remark 4.7.} Gasparis gets in fact a stronger result, that is,  a continuum of asymptotically
$l_1$ hereditarily indecomposable Banach spaces.
There is an even simpler way to obtain the result of Gasparis, if one doesn't
care for the asymptotically $l_1$ part.
In \cite{F} it was constructed a family $X_p, 1<p<+\infty$ of (uniformly convex) hereditarily
indecomposable Banach spaces with a Schauder basis. Each $X_p$ satisfies the following norm
inequality for successive vectors $x_1,\ldots,x_n$ on the basis:

$$\frac{1}{f(n)^{1/2}}(\sum_{k=1}^{n}\|x_k\|^p)^{\frac{1}{p}} \leq \|\sum_{k=1}^n x_k\| \leq (\sum_{k=1}^{n}\|x_k\|^p)^{\frac{1}{p}},$$
where as before, $f$ is defined by $f(x)=\log_2(x+1)$.
 By the same argument as in Lemma
3.7,
it follows from this inequality that for $p \neq p'$, $X_p$ and $X_{p'}$ are
totally incomparable spaces.

\vs 0.5cm

{\bf 5. Appendix.}

\vs 0.3cm

In this appendix, we give the proofs of Lemmas 3.2 and 3.3 which were
postponed in Section 3.

\vs 0.3cm

{\bf Proof of Lemma 3.2.} Since $f(x)=\hbox{log}_{2} (x+1)$ is in $\cal F$, it
follows that (i), (ii), (iii), (v) and (vi) of the definition of $\cal F$ hold
also for $f_{r}(x)= f^{r}(x)$. To show that $f_r$ belongs to $\cal F$, it only
remains to show that  the second derivate of the function  $F(x)=x/f_{r}(x)$ is negative on $[1, \infty)$.
We have,

\vs 0.1cm

 $F"(x)=\frac{r(\hbox{log}_{2}(x+1))^{r} \hbox{log}_{2} e}{(x+1)^{2}(\hbox{log}_{2}(x+1))^{2(r+1)}} \left(-x \hbox{log}_{2}(x+1)-2 \hbox{log}_{2}(x+1) +(r+1)x \hbox{log}_{2}e\right).$
\vs 0.1cm
 Consider $D(x)=-x \hbox{log}_{2}(x+1)-2 \hbox{log}_{2}(x+1) +(r+1)x
 \hbox{log}_{2}e$. Then $D(1)<0$, because $r \leq 1$ and $e^{2}<8$. Thus it is
 enough to prove that $D'(x)<0$ for all $x \geq 1$. Compute
\vs 0.1cm
 $D'(x)=\frac{1}{x+1} \left(-(x+1)\hbox{log}_{2}(x+1)-x \hbox{log}_{2}e -2 \hbox{log}_{2}e + (r+1)(x+1) \hbox{log}_{2}e\right)$.
\vs 0.1cm

 Let $H(x)=-(x+1)\hbox{log}_{2}(x+1)-x \hbox{log}_{2}e -2 \hbox{log}_{2}e + (r+1)(x+1) \hbox{log}_{2}e$.
\vs 0.1cm

 In particular, $H(1)=-2 + (2r-1) \hbox{log}{_2}e<0$ since $r \leq 1$.
Also,
 $H'(x)=- \hbox{log}_{2}(x+1) + (r-1) \hbox{log}_{2}e < 0, \ \forall x \geq
 1$, because
$r \leq 1$. Therefore $H(x)<0$ for all $x \geq 1$.  $\cqd$

\vs 0.3cm

To make clear the proof of Lemma 3.3, we stand out some  simple inequalities.

(1) $16 \hbox{log}_{2}(x^{1/2}+1) > \hbox{log}_{2} (x+1), \ \forall x \geq 1$.

(2) $(\hbox{log}_{2}(x+1))^{3/2}> \hbox{log}_{2}(x^{3}+1), \ \forall x \geq 1048576.$

(3) $\hbox{log}_{2}(x+1) \leq x^{1/4}, \ \forall x \geq 32^4.$

(4) $x^{1/4} \leq (x/2)^{1/2}, \ \forall x \geq 32^4.$

(5) $64 x^3 < e^{x/2}, \ \forall x \geq 18.$

(6) $x+1 < e^{x/2}, \ \forall x \geq 6.$

\vs 0.3cm

{\bf Proof of Lemma 3.3.} Let $J=(j_n)_{n \in \N}$ be a subset of $\N$  such that
\vspace{-0,3cm}
$$j_{1}>10^{10^{10^{400}}} \ \hbox{and} \ j_{n+1}>10^{10^{10^{4(j_n)^2}}},  \forall n \in \N.\vspace{-0,3cm}$$

 \hs 0.6cm It suffices to verify that the hypothesis of Lemma 3.1 are satisfied for this $J$ and for every $f_r$, with $r \in [1/2, 1]$. 

Fix $r \in [1/2, 1]$. By Lemma 3.2 $f_r \in \cal F$ and clearly $\hbox{log} \ \hbox{log} \ \hbox{log} \ n \geq 4m^{2}, \ \forall m, n \in J, \ m < n.$

(a) Since $f_{r}^{1/2}=f_{r/2}$, Lemma 3.2 implies that $f_{r}^{1/2} \in \cal F$.

(b) $f_{r}(j_1) > 256$ if and only if $j_{1}>2^{256^{1/r}}-1$. Since $r \geq 1/2$ this is clearly true.

(c) $exp \ exp \ j_{n} < f_{r}^{-1}(f_{r}(j_m))^{1/2}$ if and only if $\hbox{log}_{2}(1+e^{e^{j_n}})< (\hbox{log}_{2}(1+j_{m}))^{1/2}$. But this  last inequality follows from the definition of $J$.

(d) $16f_{r}(x^{1/2})>f_{r}(x)$ if and only if $16^{1/r} \hbox{log}_{2}(x^{1/2}+1)> \hbox{log}_{2}(x+1)$. Since $r\leq1$, this follows from (1).

(e) $f_{r}(p)^{3/2}>f_{r}(p^3)$ if and only if $(\hbox{log}_{2}(p+1))^{3/2}> \hbox{log}_{2}(p^{3}+1)$. This is a consequence of (2).

(f) Fix $N \in J \setminus K$ and  $x_{0}$ in the interval [log N, exp N]. The equation of the tangent t(x) to $x/{f_{r}(x)}$ at $x_0$ is
\vspace{-0,3cm}
$$t(x)= \frac{x_0}{f_{r}(x_0)} + \frac{1}{f_{r}(x_0)}(1- \frac{r x_{0}}{(x_{0}+1) \hbox{log}(x_{0}+1)})(x-x_0). \vspace{-0,3cm}$$

\vs 0.1cm

 {\it Claim 1.}
\vspace{-0,3cm}
 $$\frac{x_0}{2 (\hbox{log}_{2}(x_{0}+1))^{r+1}} \leq t(x), \ \forall x \geq 0.  \leqno (7) \vspace{-0,3cm}$$
\hs 0.6cm Indeed, since  $t(0)= (rx_{0}^{2} \ {\hbox{log}_2 \ e})/ {(x_{0}+1) (\hbox{ln}(x_{0}+1))^{r+1}}$ and $2 \leq e^{\frac{2rx_{0}}{x_{0}+1}}$, because $j_n \leq x_0$,  $\hbox{log} \ 2    \leq j_{n}/(j_{n}+1)$, for every $n \in \N$ and $r\geq 1/2$, we deduce that
\vspace{-0,3cm}
$$\frac{x_0}{2 (\hbox{log}_{2}(x_{0}+1))^{r+1}} \leq t(0). \leqno (8) \vspace{-0,3cm}$$
\hs 0.6cm Moreover, $r \leq 1$ implies that the angular coefficient  of t(x) is positive, hence $t(0) \leq t(x), \ \forall x \geq 0$.  So by (7) we conclude (8).

\vs 0.1cm

{\it Claim 2:}
 If $x < \hbox{log} \ \hbox{log} \ N$, then
\vspace{-0,3cm}
$$\frac{x}{{f_{r}(x)}^{1/2}}  \leq \frac{x_0}{2 (\hbox{log}_{2}(x_{0}+1))^{r+1}}. \leqno (9) \vspace{-0,3cm}$$
\hs 0.6cm Indeed, let $c={x_0}/{2 (\hbox{log}_{2}(x_{0}+1))^{r+1}}$. Consequently (9) holds if and only if $x^{2/r} \leq c^{2/r} \ \hbox{log}_{2}(x+1)$.

Let $d(x)=c^{2/r} \ \hbox{log}_{2}(x+1) - x^{2/r}$.

$d(1) \geq 0$ if and only if $\hbox{log}_{2}(x_{0}+1) \leq (x_{0}/2)^{1/(r+1)}$. Since $r \leq 1$, the last inequality is true because of (3) and (4).

Furthermore $d'(x)>0$ if and only if $x^{(2-r)/r} (x+1) \ 2/r <  c^{2/r} \ \hbox{log}_{2}e$, that is
\vspace{-0,3cm}
$$x^{(2-r)/r} (x+1) \ \frac{2} {r} <  \left(\frac{{x_{0}}^{1/2(r+1)}}{\hbox{log}_{2}(x_{0}+1)} \right)^{4/r}  {x_{0}}^{2/(r+1)} \ (\hbox{log}_{2}(x_{0}+1))^{2(1-r)/r} \ \frac{\hbox{log}_{2}e }{2^{2/r}}. \leqno (10) \vspace{-0,3cm}$$

Since  $1/2 \leq r \leq 1$, it follows that $(2-r)/r \leq 3$ and  $2/{r} \leq 4$. So, the first side of (10) is less than or equal to $4 x^{3} (x+1)$. On the other hand, again since $1/2 \leq r \leq 1$, we know that $1 \leq 4/r$, $1 \leq 2/(r+1)$, $0 \leq 2(1-r)/r$ and $16 \leq 1/{2^{2/r}}$. Moreover, by (7),  $1 \leq {x_{0}}^{1/2(r+1)}/{\hbox{log}_{2}(x_{0}+1)}$. Hence  the second side of (10) is greater than or equal to $x_{0}/16$.

Therefore to prove (10) it suffices to show that $64 (x+1) x^3 < x_0$.

To see this, suppose that $x < 18$, thus
\vspace{-0,3cm}
$$64 (x+1) x^3 \leq (\hbox{ln} N)^{1/2} \  (\hbox{ln} N)^{1/2} < x_{0}.  \vspace{-0,3cm}$$
\hs 0.6cm Now assume that $x \geq 18$. Hence by (5) and (6) we have
\vspace{-0,3cm}
$$64 (x+1) x^3 \leq e^{x/2} e^{x/2} < x_{0}. \vspace{-0,3cm}$$

\vs 0.1cm

 {\it Claim 3:} If $\hbox{exp} \ \hbox{exp} \ N \leq x$, then $x \geq 2x_{0}$.

\vs 0.1cm

Indeed $2x_{0} \leq 2 \ \hbox{exp} \ N \leq \hbox{exp} \ \hbox{exp} \ N \leq x$.

\vs 0.1cm

{\it Claim 4:} For every $x \geq 2x_{0}$, we have
\vspace{-0,3cm}
$$\frac{x}{4f_{r}(x_0)} \leq t(x). \vspace{-0,3cm}$$
\hs 0.6cm Indeed, consider $d(x)=t(x)- {x}/{4f_{r}(x_0)}$. Thus $d(2x_0)>0$ because  ${r x_0}/{(x_{0}+1)}< {3/2} \ \hbox{log}(x_{0}+1)$ and $d'(x)>0$ because ${r x_{0}}/{(x_{0}+1)}< {3/4} \ \hbox{log}(x_{0}+1)$.

\vs 0.1cm 

{\it Claim 5:} If $\hbox{log} \ \hbox{log} \ N \leq x$, then
\vspace{-0,3cm}
$$\frac{x}{{f_{r}(x)}^{1/2}} \leq \frac{x}{4f_{r}(x_0)}.\leqno (11) \vspace{-0,3cm}$$
\hs 0.6cm Indeed, (11) holds if and only if
\vspace{-0,3cm}
 $$16(\hbox{log}_{2}(x_{0}+1))^{2} \leq \hbox{log}_{2}(x+1) \vspace{-0,3cm}.$$
\hs 0.6cm But this is true by the definition of $J$, because
\vspace{-0,3cm}
$$16(\hbox{log}_{2}(x_{0}+1))^{2} \leq 16(\hbox{log}_{2}(e^{N}+1))^{2} \leq \hbox{log}_{2}(e^{e^{N}}+1) \leq \hbox{log}_{2}(x+1). \vspace{-0,3cm}$$
\hs 0.6cm Finally, by Claims (1), (2), (3), (4) and (5) we deduce that $x/f_{r}(x)^{1/2} \leq t(x)$, for every $x$  outside the interval $[\hbox{log} \ \hbox{log} \ N, \ \hbox{exp} \ \hbox{exp} \ N]$. $\cqd$

\
 
{\bf Acknowledgements:} We wish to thank C. Rosendal for answering many of our doubts about
descriptive set theory results.

\vs 0.1cm

 \small {Equipe d'Analyse Fonctionnelle, Universit\'e Paris 6}

 \small {Boite 186, 4, Place Jussieu, 75252, Paris Cedex 05, France}

 \small {e-mail: ferenczi@ccr.jussieu.fr}

 And

\small {Department of Mathematics, IME, University of S\~ao Paulo}

 \small {S\~ao Paulo 05315-970 Brazil}

 \small {e-mail: eloi@ime.usp.br, ferenczi@ime.usp.br}

\vs 0.6cm

\end{document}